*Research Article*

# Construction and Composition of Rooted Trees via Descent Functions

**Marco Abrate,[1] Stefano Barbero,[2] Umberto Cerruti,[2] and Nadir Murru[2]**

[1] *Dipartimento di Ingegneria Meccanica e Aerospaziale, Politecnico di Torino Corso Duca degli Abruzzi 24, 10129 Torino, Italy*
[2] *Dipartimento di Matematica, Università di Torino, Via Carlo Alberto 8, 10123 Torino, Italy*

Correspondence should be addressed to Nadir Murru; nadir.murru@gmail.com

We propose a novel approach for studying rooted trees by using functions that we will call descent functions. We provide a construction method for rooted trees that allows to study their properties through the use of descent functions. Moreover, in this way, we are able to compose rooted trees with each other. Such a new composition of rooted trees is a very powerful tool applied in this paper in order to obtain important results as the creation of new rational and Pythagorean trees.

## 1. Introduction

Trees represent one of the most important topics in the combinatorial theory. Generating trees whose nodes are elements of a structure may constitute an important step for studying the properties of the structure itself and for the structure visualization.

In this paper, some structural properties of trees will be dealt in a novel way in order to make easy the tree generation for sets satisfying specific requirements.

One of the main subjects of the present work is the so-called *descent functions*. Identifying functions of this type allows to reveal a simple and direct technique to the generation of trees on sets equipped with a weight function over a partially ordered set. Furthermore, the introduced technique yields to a substantial simplification of the study of tree levels and easily and directly finds the degree of any node.

Moreover, we introduce a *composition of trees*, which represents a completely new and very powerful tool for the generation of trees. This operation gives the feasibility of generating a new tree starting from a given pair of trees associated with a generating system once we assign a partition of the system itself. We present some examples to show the power of this tool. In particular, the composition of trees will be used to construct new examples of trees of rational numbers (rational trees) and of trees of primitive Pythagorean triples (Pythagorean trees). These applications are of special interest because nowadays only three examples of the first type [1–4] and two of the latter are known [5–7]. On the other hand, one can observe that the used techniques allow to generate an infinite number of trees of either type, and that the shown applications are merely representative.

The techniques of explicit construction of the trees described here show some features that arise from properties of the algebraic elements involved, making themselves natural and rigorous.

## 2. Tree Construction via Descent Functions

In this section, we present a new and useful approach to rooted trees by using particular functions that we will call *descent functions*. Let us start from the following.

*Definition 1.* A weight function over $X$ is a function $w : X \to Y$, for any set $X$ and $Y$, such that

(1) $Y$ is a partially ordered set, where $<$ is the relation of strict total order;

(2) $\exists \epsilon \in X : w(\epsilon) < w(x)$, for all $x \neq \epsilon \in X$;

(3) there are no infinite sequences $(x_i)_{i=0}^{\infty}$, with $x_i \in X$, such that for all $i$, $w(x_i) > w(x_{i+1})$.

It is easy to see that the element $\epsilon$ is unique. Moreover, if $Y$ satisfies the descending chain condition for sets, then condition (3) is immediately verified.

*Definition 2.* Given a weight function $w$ over $X$, a descent function is a function $\theta : X^* \to X$, where $X^* = X - \{\epsilon\}$, such that

$$w(\theta(x)) < w(x), \quad \forall x \in X^*. \tag{1}$$

In this way, for all $x \in X^*$, there exists an integer $n$ such that $\theta^n(x) = \epsilon$ (where $\theta^n = \underbrace{\theta \circ \cdots \circ \theta}_{n}$). This is sufficient to use the functions $w$ and $\theta$ for generating rooted trees.

Let us denote by $\mathcal{T}$ a rooted tree and by $\mathcal{T}_m$ its $m$th level ($\mathcal{T}_0$ contains only the root). We say that two nodes $a$ and $b$ are connected if and only if $\theta(a) = b$ or $\theta(b) = a$. If $\theta(a) = b$, then $b$ is a child node whose parent node is $a$.

*Definition 3.* Given a rooted tree $\mathcal{T}$, the level of a node $a$ is the natural number $m$ such that $a \in \mathcal{T}_m$.

The tree structure of $\mathcal{T}$ is then completely determined by the descent function $\theta$, by means of

$$\mathcal{T}_m = \{\theta(x) : \forall x \in \mathcal{T}_{m+1}\} \tag{2}$$

and vice versa one has the following:

$$\mathcal{T}_{m+1} = \{\theta^{-1}(x) : \forall x \in \mathcal{T}_m\}. \tag{3}$$

Thus, a node $a$ in $\mathcal{T}$ belongs to the $m$th level if and only if $\theta^m(a) = \epsilon$ and $\epsilon$ is the root.

*Definition 4.* We call descent tree a tree associated with a descent function $\theta$ over a set $X$, as aforementioned, having weight function $w$ over $Y$. The structure $(X, Y, w, \theta)$ is the generative tree system.

*Remark 5.* The nodes of rooted trees studied in this paper are labeled by the choice of the set $X$; that is, one and only one element of $X$ corresponds to each node of the tree. However, we will not refer in this case to the term label because in the following, we use such a term in order to classify the nodes in some classes.

*Example 6.* Let us consider the following weight and descent functions:

$$\begin{aligned} w &: \mathbb{N}_0 \longrightarrow \mathbb{N}_0, \\ w(n) &= n, \\ \theta &: \mathbb{N}_0 \setminus \{1\} \longrightarrow \mathbb{N}_0, \\ \theta(n) &= \frac{n}{2}, \quad n \text{ even}, \\ \theta(n) &= n - 1, \quad n \text{ odd}. \end{aligned} \tag{4}$$

It is straightforward to see that $w$ and $\theta$ satisfy Definitions 1 and 2, respectively, for $X = Y = \mathbb{N}$ and $\epsilon = 1$. Observing that

$$\begin{aligned} \theta^{-1}(n) &= \{n+1, 2n\}, \quad n \text{ even}, \\ \theta^{-1}(n) &= \{2n\}, \quad n \text{ odd}, \end{aligned} \tag{5}$$

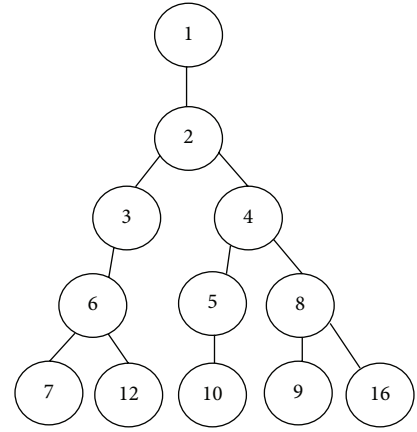

Figure 1

starting from the root 1, it is easy to construct the tree, represented up to the 4th level in Figure 1.

*Example 7.* Let $X$ be the set of vectors $(a_1, \ldots, a_t)$ such that $t \geq 1$, $a_i \in \mathbb{N}_0$ for every $i$, and $a_1 \geq a_2 \geq \cdots \geq a_t$. Moreover, let us consider $Y = \mathbb{N}$ and the weight function $w : X \to Y$ defined by

$$w(a_1, \ldots, a_t) = t + \sum_{i=1}^{t} a_i. \tag{6}$$

And thus $\epsilon = (1)$. Let us use the following descent function:

$$\begin{aligned} \theta(a_1, \ldots, a_t) &= (a_1, \ldots, a_{t-1}), \quad \text{if } a_t = 1, \\ \theta(a_1, \ldots, a_t) &= (a_1, \ldots, a_{t-1}, a_t - 1), \quad \text{if } a_t > 1, \end{aligned} \tag{7}$$

whose inverse is

$$\begin{aligned} \theta^{-1}(a_1, \ldots, a_t) &= (a_1, \ldots, a_t, 1), \quad \text{if } a_{t-1} = a_t = 1, \\ \theta^{-1}(a_1, \ldots, a_t) &= \{(a_1, \ldots, a_{t-1}, a_t + 1), (a_1, \ldots, a_t, 1)\}, \\ &\quad \text{otherwise}. \end{aligned} \tag{8}$$

The resulting tree is

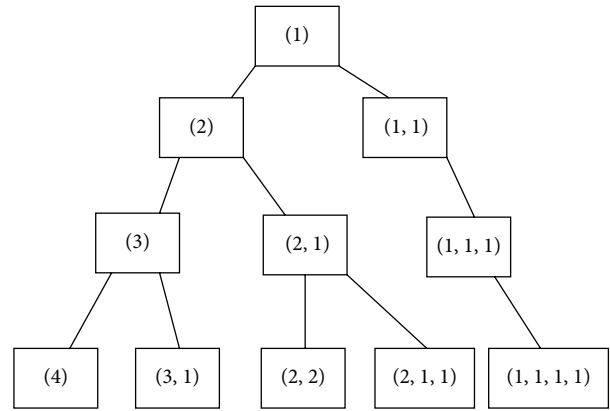

where the level $\mathcal{T}_m$ contains all the partitions of the integer $m$.

*Definition 8.* We call natural descent function of a rooted tree $\mathcal{T}$ the function $\nu_\mathcal{T}$ such that $\nu_\mathcal{T}(a)$ is the parent node of $a$, for any node $a \in \mathcal{T}$.

*Definition 9.* Let us define the set $U = \{n \in \mathbb{N} : \text{if } p|n \text{ and } q < p \text{ then } q|n\}$.

*Definition 10.* Let $\mathcal{U}$ be a descent tree, whose nodes are labeled by $U$. The structure of $\mathcal{U}$ is determined by the following descent function:

$$\theta_U(a) = \theta\left(\prod_{i=1}^m p_i^{h_i}\right) = \frac{a}{p_m^{h_m}}. \qquad (9)$$

The tree $\mathcal{U}$ is an infinite tree, whose nodes have infinite degree.

**Theorem 11.** *Given any countable rooted tree $\mathcal{T}$, then there exists an injective function $j : \mathcal{T} \to \mathcal{U}$, such that the following diagram is commutative:*

$$\begin{array}{ccc} \mathcal{T} & \xrightarrow{\nu_\mathcal{T}} & \mathcal{T} \\ \downarrow j & & \downarrow j \\ \mathcal{U} & \xrightarrow{\theta_U} & \mathcal{U} \end{array} \qquad (10)$$

The function $j$ of the theorem can be defined recursively. The idea is to assign an element of $x \in R$ to a node of $\mathcal{T}$ that belongs to the $k$th floor, such that the prime $p_k$ is the greatest divisor of $x$. In particular, if $\epsilon$ is the root of $\mathcal{T}$, then $j(\epsilon) = 1$. Moreover, let us suppose that $j$ is defined over all elements of $\mathcal{T}$ until the floor $\mathcal{T}_k$. A node $a \in \mathcal{T}_{k+1}$ is the child of a node $b \in \mathcal{T}_k$. If $j(b) = p_1^{e_1} p_2^{e_2} \cdots p_k^{e_k}$ and $a$ is the $s$th child of $b$, then

$$j(a) = p_1^{e_1} p_2^{e_2} \cdots p_k^{e_k} p_{k+1}^s. \qquad (11)$$

In this way, any countable rooted tree is contained in the *universal* tree $\mathcal{U}$. Indeed, given a countable rooted tree $\mathcal{T}$, it is determined by its natural descent function $\nu_\mathcal{T}$. Furthermore, $\nu_\mathcal{T}$ is the restriction of $\theta_U$ to $\mathcal{T}$.

*Example 12.* Let us consider the following tree without label:

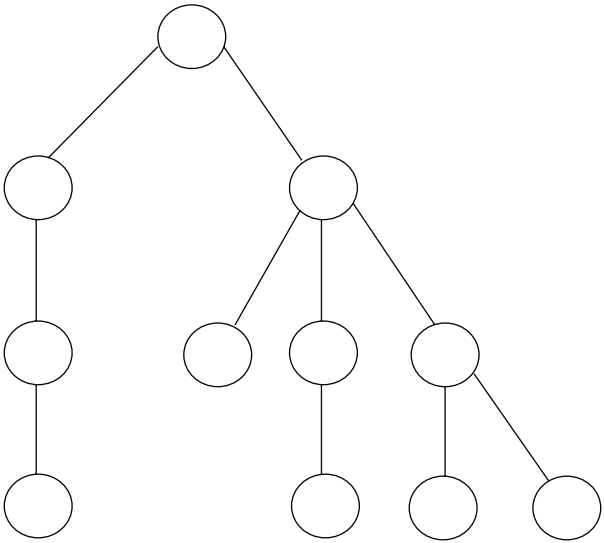

We can assign labels, embedding the tree into $\mathcal{U}$ in the following way:

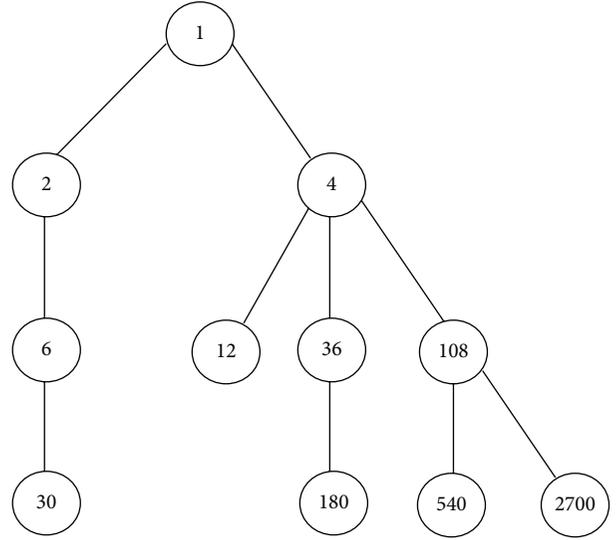

that is,

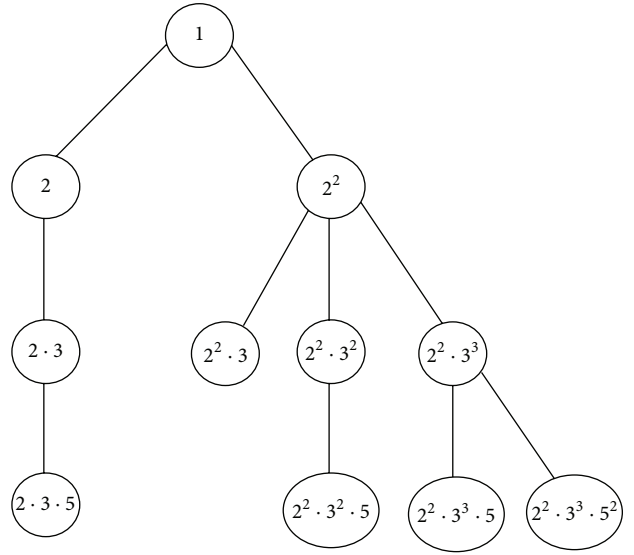

The structure of the tree is determined by the descent function (9). For example, we can determine the paths of the node 30 as follows:

$$\theta_U(30) = \theta_U(2 \cdot 3 \cdot 5) = 2 \cdot 3 = 6,$$
$$\theta_U(6) = \theta_U(2 \cdot 3) = 2, \qquad \theta_U(2) = 1. \qquad (12)$$

*Example 13.* The tree given in Example 6 can be labeled as follows:

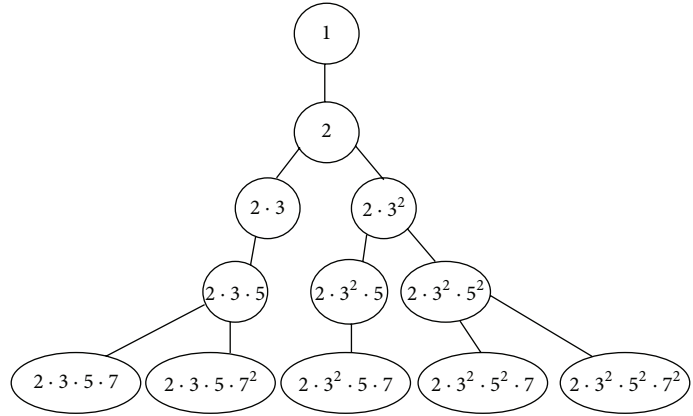

And consequently it can be derived by using the descent function (9). We can derive such a tree starting from $U$ and $\theta_U$. Given $n \in U$, let $u(n)$ and $h(n)$ be the exponent of the greatest prime that divides $n$ and $h(n)$ the index of the greatest prime that divides $n$ in the sequence of prime numbers, respectively. Thus, the previous tree is defined by the following *laws of growth*:

(1) 1 generates 2;

(2) if $h(n) = k$ and $u(n) = 1$, then $n$ generates the children $np_{k+1}$ and $np_{k+1}^2$;

(3) if $h(n) = k$ and $u(n) = 2$, then $n$ generates the child $np_{k+1}$.

Now we focus on a class of rooted trees which satisfy the following.

*Definition 14.* Let us consider the rooted tree $\mathcal{T}$, whose nodes are labeled by the set $X = \bigsqcup_{i=1}^{v} A_i$ (disjoint union). We call types the sets $A_i$. The fact that $a \in A_i$, that is, $a$ is a node of the type $A_i$, determines the number of its children and their type. We call $\mathcal{T}$ a typed tree, and if $a \in A_i$, we indicate with $s_{ij}$ the number of its children that are in $A_j$ (i.e., the number of its children that have type $A_j$).

In a typed tree, any node belongs to a type; that is, the number of children and their type is known. A similar class of tree is studied in [8, 9], where the nodes of the tree are labeled considering the number of children and their label.

*Remark 15.* In the following, we directly use typed trees to derive the order of the levels as recurrent sequences. A sequence $(a_n)_{n=0}^{\infty}$ recurs with characteristic polynomial

$$f(t) = t^m - \sum_{i=1}^{m} h_i t^{m-i}, \tag{13}$$

if $(a_n)_{n=0}^{\infty}$ satisfies the linear homogeneous recurrent relation

$$a_n = \sum_{i=1}^{m} h_i a_{n-i}, \quad \forall n \geq m, \tag{14}$$

given the initial conditions $a_0, \ldots, a_{m-1}$. We also recall that the ordinary generating function related to a linear recurrent sequence is $F(t) = u(t)/f^R(t)$, where $f^R(t) = 1 - \sum_{i=1}^{m} h_i t^i$ is the reflected polynomial of $f(t)$ and the polynomial $u(t) = \sum_{i=1}^{m-1} u_i t^i$ is determined from initial conditions (see, e.g., [10] for a survey on recurrent sequences).

Using previous notation, we can define the *matrix of types* as $G = (s_{ij})$. Such a matrix easily allows to determine the counting sequence of the levels as a linear recurrent sequence.

**Proposition 16.** *The sequence $(|\mathcal{T}_m|)_{m=0}^{\infty}$ is a linear recurrent sequence whose characteristic polynomial is the characteristic polynomial of $G$.*

*Proof.* Let us define $u_i^{(m)} = |\{x \in \mathcal{T}_m : x \in A_i\}|$; then $|\mathcal{T}_m| = \sum_{i=1}^{v} u_i^{(m)}$. Moreover, defining the distribution of the elements of $\mathcal{T}_m$ as

$$\pi(m) = \left(u_1^{(m)}, \ldots, u_v^{(m)}\right), \tag{15}$$

it is immediate that

$$\pi(m) G = \pi(m+1). \tag{16}$$

It follows that any sequence $(u_i^{(m)})_{m=0}^{\infty}$ recurs with the characteristic polynomial of $G$. Since $|\mathcal{T}_m|$ is the sum of the $u_i^{(m)}$'s, then it is wellknown that $|\mathcal{T}_m|$ is a linear recurrent sequence with the same characteristic polynomial of the $u_i^{(m)}$'s. □

*Example 17.* Let us consider the tree constructed in Example 6. We consider the partition $\mathbb{N} = A_1 \cup A_2$ where $A_1 = \{n \in \mathbb{N} : n \text{ odd}\}$ and $A_2 = \{n \in \mathbb{N} : n \text{ even}\}$. Then $G = \begin{pmatrix} 0 & 1 \\ 1 & 1 \end{pmatrix}$: odd numbers have only one child which is an even number; even numbers have two children (one odd and one even). As a consequence, the sequence of the levels orders of such tree is the Fibonacci sequence. Indeed, the characteristic polynomial of $G$ is $x^2 - x - 1$, and it is straightforward to check that $|\mathcal{T}_0| = 1$, $|\mathcal{T}_1| = 1$.

*Example 18.* Let us consider the tree constructed in Example 7. This is an example of tree constructed via descent function which is not a typed tree. By Proposition 16 it follows that in a typed tree the orders of the levels are a linear recurrent sequence. In Example 7, $|\mathcal{T}_m|$ is the number of partitions of the integer $m$, and it is well known that such a sequence is not recurrent. Indeed, it is wellknown that the generating function of the sequence counting the number of partition of the integers is

$$\prod_{k>0} \frac{1}{1-x^k}. \qquad (17)$$

That is not a rational function (see, e.g., sequence A000041 in [11]).

However, it is interesting to see that only changing the set $X$, the descent tree alters its structure and properties. Let us use the same set $Y$ and functions $w, \theta$ of the Example 7, but we change $X$, taking the set of the vectors $(a_1, \ldots, a_t)$ such that $t \geq 1$ and $a_i \in \mathbb{N}$ for every $i$. In this case, the inverse descent function is

$$\theta^{-1}(a_1, \ldots, a_t) = \{(a_1, \ldots, a_t, 1), (a_1, \ldots, a_{t-1}, a_t + 1)\}, \qquad (18)$$

and the descent tree is

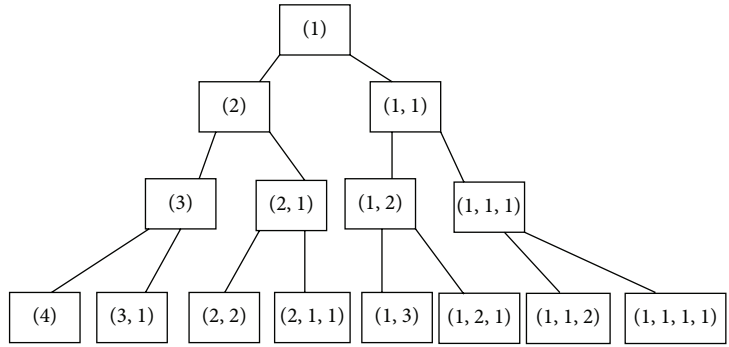

Setting $X = A_1 \cup A_2$, where $A_1 = \{a \in X : \text{last component} = 1\}$ and $A_2 = \{a \in X : \text{last component} > 1\}$, this is a generating tree. We easily see that $G = \begin{pmatrix} 1 & 1 \\ 1 & 1 \end{pmatrix}$, whose characteristic polynomial is $x^2 - 2x$. Thus, we have $|\mathcal{T}_m| = 2^m$, for any $m \geq 0$.

An important consequence of the use of descent functions for constructing trees is the possibility of composing trees to each other in order to obtain new ones. The composition of two descent trees is simple, and it is highlighted in the following.

*Definition 19.* Let $\mathcal{T}, S$ be descent trees created by the generative tree systems $(X, Y, w, f)$ and $(X, Y, w, g)$, respectively. Let us consider any partition $X = A \cup B$ formed by two subsets. We define $\mathcal{R} = \mathcal{T} \circ_{A,B} S$ as the composed descent tree whose associated generative system is $(X, Y, w, j_{A,B})$, where

$$j_{A,B} : X^* \longrightarrow X,$$
$$j_{A,B}(x) = f(x), \quad \forall x \in A, \qquad (19)$$
$$j_{A,B}(x) = g(x), \quad \forall x \in B.$$

It is easy to see that since $f$ and $g$ are descent functions, the function $j_{A,B}$ is still a descent function, and therefore $\mathcal{R}$ is still a descent tree. By Definition 19, it is clear that tree composition depends on the partition of $X$ and on the order of the subsets $A$ and $B$ considered in (19). This means that any partition generates two different trees; that is,

$$\mathcal{T} \circ_{A,B} S \neq \mathcal{T} \circ_{B,A} S. \qquad (20)$$

In the next sections, we apply these tools to rational trees and Pythagorean trees.

## 3. Applications and Examples

*3.1. Rational Trees.* A rational tree is a tree such that $X = \mathbb{Q}^+$, where every positive rational number appears in the tree exactly once and each node generates children following some explicit rule. Actually, only three rational trees are known: the Kepler tree [3], the Calkin-Wilf tree [2], and the Stern-Brocot tree [1, 4].

The oldest, but the less famous, is the Kepler tree, introduced by Kepler in Harmonices Mundi in 1619. The Kepler tree is a binary tree, that is; each node has two children by means of the following rule:

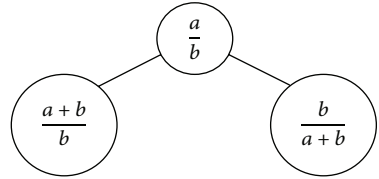

We can derive the Kepler tree via descent functions, by setting

$$X = \mathbb{Q}^+, \qquad Y = \mathbb{N},$$
$$w\left(\frac{a}{b}\right) = a + b, \qquad \epsilon = \frac{1}{1}, \qquad (21)$$
$$\theta_K\left(\frac{a}{b}\right) = \begin{cases} \dfrac{a-b}{b}, & a > b, \\ \dfrac{b-a}{a}, & a < b. \end{cases}$$

*Remark 20.* The fractions $a/b$ in rational trees are considered with $a$, $b$ coprime.

It is easy to see that $\theta_K$ generates the Kepler tree:

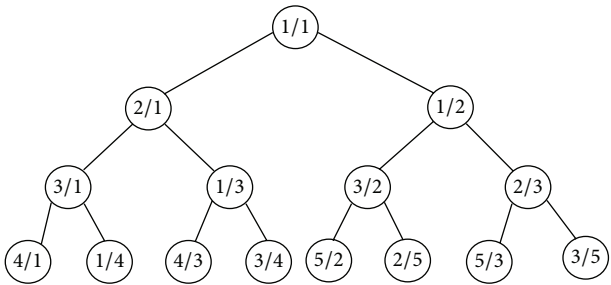

A similar tree, recently studied, is the Calkin-Wilf tree where the construction rule is given by

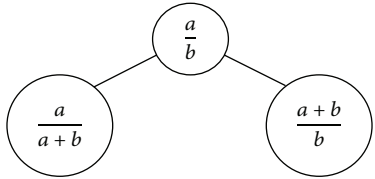

We can derive the corresponding generative system tree by means of

$$X = \mathbb{Q}^+, \qquad Y = \mathbb{N},$$
$$w\left(\frac{a}{b}\right) = a + b, \qquad \epsilon = \frac{1}{1}, \qquad (22)$$
$$\theta_C\left(\frac{a}{b}\right) = \begin{cases} \dfrac{a-b}{b}, & a > b, \\ \dfrac{a}{b-a}, & a < b, \end{cases}$$

from which we have the Calkin-Wilf tree:

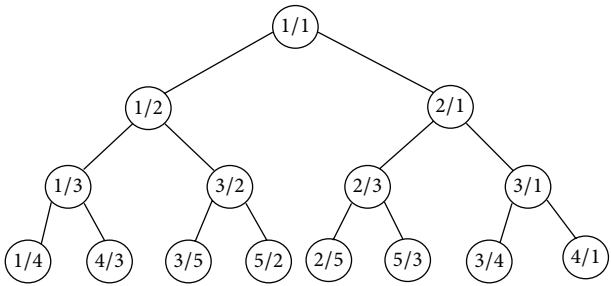

Now, we can apply the composition of descent functions stated in Definition 19 to the generative systems tree $(X, Y, w, \theta_K)$ and $(X, Y, w, \theta_C)$ obtaining a new generative system tree $(X, Y, w, \theta)$, where, for example, we can construct $\theta$ as follows:

$$\theta\left(\frac{a}{b}\right) = \begin{cases} \theta_K\left(\dfrac{a}{b}\right), & b \text{ even}, \\ \theta_C\left(\dfrac{a}{b}\right), & b \text{ odd}, \end{cases}$$
$$= \begin{cases} \dfrac{a-b}{b}, & a > b, \\ \dfrac{b-a}{a}, & a < b \text{ and } b \text{ even}, \\ \dfrac{a}{b-a}, & a < b \text{ and } b \text{ odd}, \end{cases} \qquad (23)$$

from which we can derive its inverse in order to construct the tree starting from the root $1/1$ as follows:

$$\theta^{-1}\left(\frac{a}{b}\right) = \left\{\frac{a}{a+b}, \frac{a+b}{b}\right\}, \quad \text{if } a + b \text{ odd},$$
$$\theta^{-1}\left(\frac{a}{b}\right) = \left\{\frac{a+b}{b}, \frac{b}{a+b}\right\}, \quad \text{if } a + b \text{ even}. \qquad (24)$$

The generative system tree $(X, Y, w, \theta)$ generates the following tree:

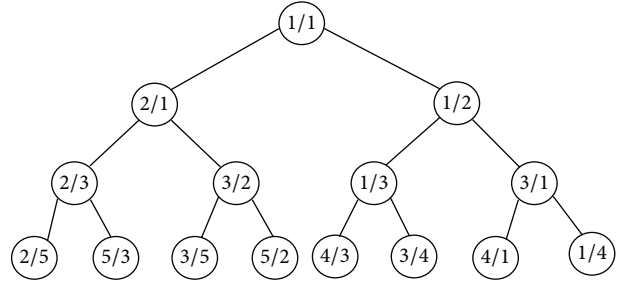

Such a tree is different from the Kepler, Calkin-Wilf, and Stern-Brocot trees. Moreover, every positive rational number appears exactly once in the tree. Indeed, given any rational number $a/b$, there exists a positive integer $n$ such that $\theta^n(a/b) = \epsilon = 1/1$, since $\theta$ is a descent function for the system $(X, Y, w, \theta)$ by construction.

**Theorem 21.** *The generative system tree $(X, Y, w, \theta)$, where*

$$X = \mathbb{Q}^+, \qquad Y = \mathbb{N},$$
$$w\left(\frac{a}{b}\right) = a + b, \quad \forall \frac{a}{b} \in X,$$
$$\theta\left(\frac{a}{b}\right) = \begin{cases} \dfrac{a-b}{b}, & a > b, \\ \dfrac{b-a}{a}, & a < b \text{ and } b \text{ even}, \\ \dfrac{a}{b-a}, & a < b \text{ and } b \text{ odd}, \end{cases} \qquad (25)$$

*generates a rational tree.*

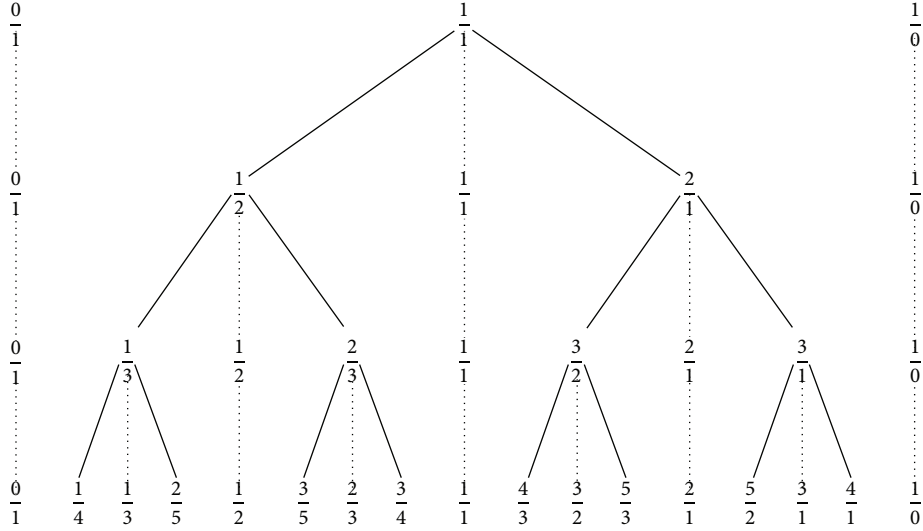

Figure 2

*Remark 22.* Using different partitions of $\mathbb{Q}^+$ into two subsets, we can apply Definition 19 in order to compose Kepler and Calkin-Wilf trees to obtain a new rational tree. This is a very important result because in this way, we are able to construct infinite rational trees.

Finally, we examine the Stern-Brocot tree from the point of view of the descent functions. The Stern-Brocot tree is a rational tree made up of mediants $(a + b)/(c + d)$ of two fractions $a/c$ and $b/d$, as described in Figure 2. It is not immediate to obtain the corresponding generative system tree (despite the Kepler and Calkin-Wilf cases). However, if we consider the $2 \times 2$ invertible matrices with determinant 1, we can find a descent function from which we can derive the Stern-Brocot tree. Precisely, we take the following:

$$X = SL(2, \mathbb{N}), \quad Y = \mathbb{N},$$
$$w\left(\begin{bmatrix} a & b \\ c & d \end{bmatrix}\right) = a+b+c+d, \quad \epsilon = \begin{bmatrix} 1 & 0 \\ 0 & 1 \end{bmatrix}, \quad (26)$$

$$\theta_S\left(\begin{bmatrix} a & b \\ c & d \end{bmatrix}\right) = \begin{cases} \begin{bmatrix} a-b & b \\ c-d & d \end{bmatrix}, & \begin{bmatrix} a \\ c \end{bmatrix} > \begin{bmatrix} b \\ d \end{bmatrix} \\ \begin{bmatrix} a & b-a \\ c & d-c \end{bmatrix}, & \begin{bmatrix} a \\ c \end{bmatrix} < \begin{bmatrix} b \\ d \end{bmatrix}, \end{cases} \quad (27)$$

where $\begin{bmatrix} a \\ c \end{bmatrix} > \begin{bmatrix} b \\ d \end{bmatrix}$ means

$$a \geq b, \quad c \geq d. \quad (28)$$

*Remark 23.* Observe that, since the matrices considered ealier are invertible, equalities in (28) cannot hold simultaneously. Moreover, given a $\begin{bmatrix} a & b \\ c & d \end{bmatrix} \in SL(2, \mathbb{N})$, $\begin{bmatrix} a \\ c \end{bmatrix} > \begin{bmatrix} b \\ d \end{bmatrix}$ or $\begin{bmatrix} a \\ c \end{bmatrix} < \begin{bmatrix} b \\ d \end{bmatrix}$ holds, except for $\epsilon$.

It is easy to observe that

$$\theta_S^{-1}\left(\begin{bmatrix} a & b \\ c & d \end{bmatrix}\right) = \left\{\begin{bmatrix} a+b & b \\ c+d & d \end{bmatrix}, \begin{bmatrix} a & a+b \\ c & c+d \end{bmatrix}\right\}, \quad (29)$$

using the following transformation,

$$\begin{bmatrix} a & b \\ c & d \end{bmatrix} \longrightarrow \frac{a+b}{c+d}. \quad (30)$$

We can construct the Stern-Brocot tree by using the generative system tree $(X, Y, w, \theta_S)$.

*Remark 24.* The use of the descent function $\theta_S$ provides a new explicit representation of the Stern-Brocot tree.

In [12], a method involving continued fractions has been proposed in order to locate the position of any rational number over the Stern-Brocot tree. Here, we highlight an iterative method (involving the generative system tree $(X, Y, w, \theta_S)$) to find the path of any rational number over the Stern-Brocot tree.

Given any rational number $e/f$ in the Stern-Brocot tree (with $e \neq 1$), then it is the mediant of the fractions $a/c$, $b/d$, where

$$a = f^{-1} \bmod e, \quad b = e - a,$$
$$c = \frac{af - 1}{e}, \quad d = f - c. \quad (31)$$

Indeed, if $e/f$ is the mediant of $a/c$ and $b/d$ we must have

$$a + b = e, \quad c + d = f, \quad ad - bc = 1, \quad (32)$$

and we obtain

$$af - 1 = ac + ad - ad + bc = ac + bc = ce, \quad (33)$$

from which

$$b = e - a, \qquad c = \frac{af - 1}{e}, \qquad d = f - c. \qquad (34)$$

We only need to find $a$. Since $e$ and $f$ are coprime, $f \in \mathbb{Z}_e^*$ and we obtain

$$a = f^{-1} \bmod e. \qquad (35)$$

If the rational number is $1/f$ the path on the Stern-Brocot tree is trivial to find.

Thus, given any rational number $e/f$, we are able to construct the corresponding matrix $\begin{bmatrix} a & b \\ c & d \end{bmatrix} \in SL(2, \mathbb{N})$ such that $a + b = e$ and $c + d = f$. Using the descent function $\theta_S$, we have

$$\theta_S\left(\begin{bmatrix} a & b \\ c & d \end{bmatrix}\right) = \begin{bmatrix} a' & b' \\ c' & d' \end{bmatrix}, \qquad (36)$$

according to (27), where $(a'+b')/(c'+d')$ is the parent node of $e/f$ in the Stern-Brocot tree. The iteration of the descent function $\theta_S$ provides the path of the rational number $e/f$.

*Example 25.* Let us consider the rational number $11/8$; then by (31) we have that in the Stern-Brocot tree $11/8$ is the mediant of

$$\frac{7}{5}, \quad \frac{4}{3}, \qquad (37)$$

$$\theta_S\left(\begin{bmatrix} 7 & 4 \\ 5 & 3 \end{bmatrix}\right) = \begin{bmatrix} 3 & 4 \\ 2 & 3 \end{bmatrix}. \qquad (38)$$

Consequently the parent node of $11/8$ is $7/5$. Now, we know from (38) that $7/5$ is the mediant of $3/2$ and $4/3$, and

$$\theta_S\left(\begin{bmatrix} 3 & 4 \\ 2 & 3 \end{bmatrix}\right) = \begin{bmatrix} 3 & 1 \\ 2 & 1 \end{bmatrix}, \qquad (39)$$

and the parent node of $7/5$ is $4/3$. Again, we have

$$\theta_S\left(\begin{bmatrix} 3 & 1 \\ 2 & 1 \end{bmatrix}\right) = \begin{bmatrix} 2 & 1 \\ 1 & 1 \end{bmatrix}. \qquad (40)$$

The parent node of $4/3$ is $3/2$. Similarly, we obtain that the parent node of $3/2$ is $2/1$, and finally we arrive to the root $1/1$. Thus the path from the rational number $11/8$ to the root in the Stern-Brocot tree is

$$\left(\frac{11}{8}, \frac{7}{5}, \frac{4}{3}, \frac{3}{2}, \frac{2}{1}, \frac{1}{1}\right). \qquad (41)$$

*3.2. Pythagorean Trees.* Pythagorean trees are trees of Pythagorean triples; that is, each node is a triple $[a, b, c]$ of positive integers that satisfy the equation $a^2 + b^2 = c^2$. In the following, we describe the two known Pythagorean trees in terms of descent functions introduced in Section 2. Moreover, we will find two new Pythagorean trees by simply applying the tree composition presented in Definition 19.

We consider the sets $Y = 2\mathbb{N}$ and $X$ of all the couples $(a, b)$ where $a > b \geq 1$ are odd coprime positive integers. The function

$$\begin{aligned} w : X &\longrightarrow Y, \\ w((a,b)) &= a + b \end{aligned} \qquad (42)$$

satisfies the Definition 1 of weight function. As a matter of fact, $Y$ is a partially ordered set with respect to the relation of strict total order $<$ between integers, and we have $\epsilon = (3, 1)$ such that $w(\epsilon) < w((a,b))$ for all $(a,b) \neq \epsilon \in X$, and also the third condition of Definition 1 is clearly satisfied. The two systems $(X, Y, w, \theta_1)$ and $(X, Y, w, \theta_2)$ where

$$\theta_1((a,b)) = \begin{cases} (a - 2b, b) & \text{if } a > 3b, \\ (b, a - 2b) & \text{if } 2b < a < 3b, \\ (b, 2b - a) & \text{if } b < a < 2b, \end{cases}$$

$$\theta_2((a,b))$$
$$= \begin{cases} \left(\dfrac{a+b}{2}, b\right) & \text{if } \dfrac{a+b}{2} \text{ is odd}, \\ \left(\max\left(\dfrac{a-b}{2}, b\right), \min\left(\dfrac{a-b}{2}, b\right)\right) & \text{if } \dfrac{a-b}{2} \text{ is odd}, \end{cases}$$
$$(43)$$

are generative tree systems. Indeed, from (43), we immediately have for all $x \in X^*$ that $w(\theta_1(x)) < w(x)$ and $w(\theta_2(x)) < w(x)$, so the functions $\theta_1$ and $\theta_2$ are descent functions. Moreover, a simple calculation shows that

$$\theta_1^{-1}((a,b)) = \{(a + 2b, b), (2a + b, a), (2a - b, a)\}, \qquad (44)$$

$$\theta_2^{-1}((a,b)) = \{(2a - b, b), (2a + b, b), (a + 2b, a)\}. \qquad (45)$$

So we have a complete description of the two trees $\mathscr{PT}^1$ and $\mathscr{PT}^2$ corresponding to the generative tree systems $(X, Y, w, \theta_1)$ and $(X, Y, w, \theta_2)$, respectively. We show $\mathscr{PT}^1$ to the third level as follows:

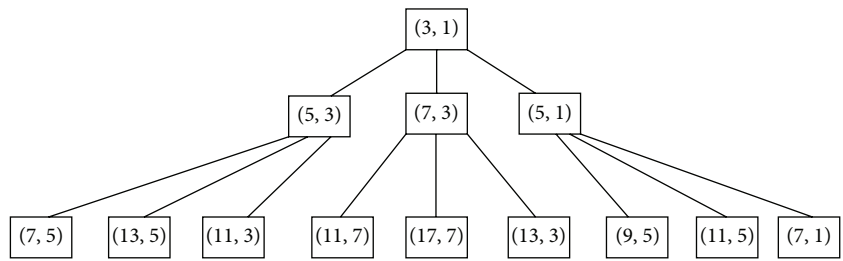

and $\mathscr{PT}^2$ to the third level as follows:

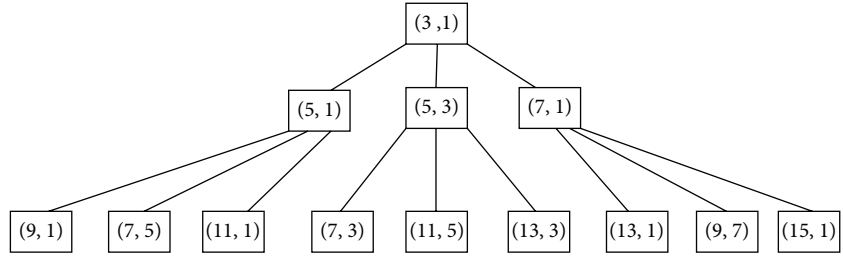

These trees are strictly related to the Pythagorean triples. As a matter of fact, we can observe their correspondence with the two trees presented in [7] where the second columns of the matrices in the nodes are exactly the nodes of our trees. Here every node $(a, b)$ represents the primitive Pythagorean triple $(ab, (a^2 - b^2)/2, (a^2 + b^2)/2)$. The tree $\mathscr{PT}^1$ is the Barning-Hall tree [5, 6], while $\mathscr{PT}^2$ is the new Pythagorean tree presented by Price [7]. Combining $\mathscr{PT}^1$ and $\mathscr{PT}^2$ by means of the composition rule stated in Definition 19, we can generate two new Pythagorean trees. As a first step, we have to find two disjoint subsets in $X$; for example, we can consider the most natural partition of the set of even numbers; namely,

$$\begin{aligned} A_1 &= \{(a,b) \in X : a = b \bmod 4\} \\ &= \{(a,b) \in X : a + b = 2 \bmod 4\}, \\ A_2 &= \{(a,b) \in X : a \neq b \bmod 4\} \\ &= \{(a,b) \in X : a + b = 0 \bmod 4\}. \end{aligned} \quad (46)$$

Thanks to this partition of $X$, we obtain the two compositions

$$\Theta_1((a,b)) = \begin{cases} \theta_1((a,b)) & \text{if } (a,b) \in A_1, \\ \theta_2((a,b)) & \text{if } (a,b) \in A_2, \end{cases}$$

$$\Theta_2((a,b)) = \begin{cases} \theta_1((a,b)) & \text{if } (a,b) \in A_2, \\ \theta_2((a,b)) & \text{if } (a,b) \in A_1, \end{cases} \quad (47)$$

which clearly are descent functions. We easily determine their inverses using the definitions of $\theta_1^{-1}$ and $\theta_2^{-1}$ as follows. When $(a, b) \in A_1$, we have for the elements belonging to the set (44)

$$w((a + 2b, b)) = a + 3b = 0 \bmod 4 \implies (a + 2b, b) \in A_2,$$
$$w((2a + b, a)) = 3a + b = 0 \bmod 4 \implies (2a + b, a) \in A_2,$$
$$w((2a - b, a)) = 3a - b = 2 \bmod 4 \implies (2a - b, a) \in A_1 \quad (48)$$

for the elements of the set (45),

$$w((2a - b, b)) = 2a = 2 \bmod 4 \implies (2a - b, b) \in A_1,$$
$$w((2a + b, b)) = 2a + 2b = 0,$$
$$\bmod 4 \implies (2a + b, b) \in A_2, \quad (49)$$
$$w((a + 2b, a)) = 2a + 2b = 0,$$
$$\bmod 4 \implies (a + 2b, a) \in A_2.$$

On the other hand, when $(a, b) \in A_2$ we have for the elements belonging to the set (44)

$$w((a + 2b, b)) = a + 3b = 2 \bmod 4 \implies (a + 2b, b) \in A_1,$$
$$w((2a + b, a)) = 3a + b = 2 \bmod 4 \implies (2a + b, a) \in A_1,$$
$$w((2a - b, a)) = 3a - b = 0 \bmod 4 \implies (2a - b, a) \in A_2; \quad (50)$$

for the elements of the set (45), the situation is the same as in (49). Therefore, if we use $\Theta_1$, we reach $(a, b) \in A_1$ from the node $(2a-b, a) \in A_1$ using $\theta_1$ or from the nodes $(2a+b, b), (a+2b, a) \in A_2$ using $\theta_2$. Moreover, we reach $(a, b) \in A_2$ from the nodes $(a+2b, b), (2a+b, b) \in A_1$ using $\theta_1$, and from the nodes $(2a + b, b), (a + 2b, a) \in A_2$ using $\theta_2$. This allows us to invert $\Theta_1$ in this way as follows:

$$\Theta_1^{-1}((a,b))$$
$$= \begin{cases} \{(2a - b, a), (2a + b, b), (a + 2b, a)\} \\ \qquad\qquad\qquad\qquad\qquad \text{if } (a,b) \in A_1, \\ \{(a + 2b, b), (2a + b, a), (2a + b, b), (a + 2b, a)\} \\ \qquad\qquad\qquad\qquad\qquad \text{if } (a,b) \in A_2. \end{cases} \quad (51)$$

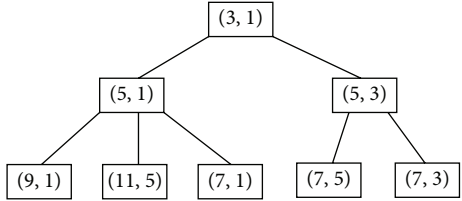

Figure 3

With similar considerations on $\Theta_2$, we obtain

$$\Theta_2^{-1}((a,b)) = \begin{cases} \{(2a-b,b),(2a+b,a),(a+2b,b)\} & \text{if } (a,b) \in A_1, \\ \{(2a-b,a),(2a-b,b)\} & \text{if } (a,b) \in A_2. \end{cases} \quad (52)$$

The system $(X, Y, w, \Theta_1)$ generates the tree $\mathscr{PT}^3 = \mathscr{PT}^1 \circ_{A_1, A_2} \mathscr{PT}^2$ which we show to the third level

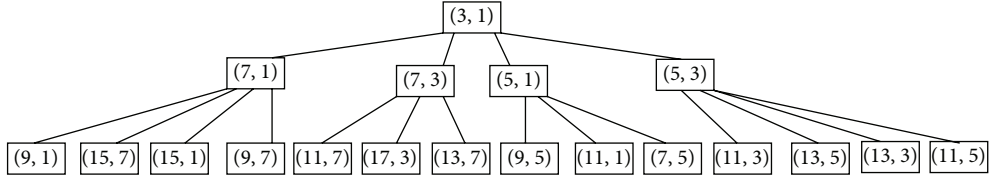

The system $(X, Y, w, \Theta_2)$ gives rise to the tree $\mathscr{PT}^4 = \mathscr{PT}^1 \circ_{A_2, A_1} \mathscr{PT}^2$ whose first three levels appear as in Figure 3.

The partition of $X$ in the two subsets $A_1$ and $A_2$ allows us to complete our analysis of the four Pythagorean trees that we have considered.

*Remark 26.* We point out some interesting results which are straightforward applications of our approach. For each of the four Pythagorean trees described ealier, we are able to determine the related generating matrix and consequently to study the sequence of the order of the $m$th level for $m \geq 0$. In fact, from relations (48),(49), and (50) it easily follows that

(i) the generating matrix for $\mathscr{PT}^1$ is $G_1 = \begin{pmatrix} 1 & 2 \\ 2 & 1 \end{pmatrix}$ whose characteristic polynomial is $x^2 - 4x + 3$, with roots $x = 3$ and $x = 1$;

(ii) the generating matrix for $\mathscr{PT}^2$ is $G_2 = \begin{pmatrix} 1 & 2 \\ 1 & 2 \end{pmatrix}$ whose characteristic polynomial is $x^2 - 3x$ with roots $x = 0$ and $x = 3$;

(iii) the generating matrix for $\mathscr{PT}^3$ is $G_3 = \begin{pmatrix} 1 & 2 \\ 2 & 2 \end{pmatrix}$ whose characteristic polynomial is $x^2 - 3x - 2$ with roots $x = (3 + \sqrt{17})/2$ and $x = (3 - \sqrt{17})/2$;

(iv) the generating matrix for $\mathscr{PT}^4$ is $G_4 = \begin{pmatrix} 1 & 2 \\ 1 & 1 \end{pmatrix}$ whose characteristic polynomial is $x^2 - 2x - 1$ with roots $x = 1 + \sqrt{2}$ and $x = 1 - \sqrt{2}$.

Clearly, the $m$th level order of $\mathscr{PT}^1$ and $\mathscr{PT}^2$ is the same as follows:

$$\left|\mathscr{PT}_m^1\right| = \left|\mathscr{PT}_m^2\right| = 3^m. \quad (53)$$

The $m$th level order of $\mathscr{PT}^3$ is

$$\left|\mathscr{PT}_m^3\right| = \frac{5+\sqrt{17}}{2\sqrt{17}}\left(\frac{3+\sqrt{17}}{2}\right)^m - \frac{5-\sqrt{17}}{2\sqrt{17}}\left(\frac{3-\sqrt{17}}{2}\right)^m. \quad (54)$$

We observe that the sequence of level orders has ordinary generating function (see [10]) given by

$$F_{\mathscr{PT}^3}(t) = \frac{1+t}{1-3t-2t^2} = 1 + 4t + 14t^2 + 50t^3 + 178t^4 \\ + 634t^5 + 2258t^6 + 8042t^7 \\ + 28642t^8 + O(t^9), \quad (55)$$

and $(|\mathscr{PT}_m^3|)_{m=0}^{+\infty} = (1, 4, 14, 50, 178, 634, 2258, 8042, 28642, \ldots)$ is a sequence not present in the OEIS.

Finally, the $m$th level order of $\mathscr{PT}^4$ is

$$\left|\mathscr{PT}_m^4\right| = \frac{1}{2\sqrt{2}}\left[\left(1+\sqrt{2}\right)^{m+1} - \left(1-\sqrt{2}\right)^{m+1}\right]. \quad (56)$$

The sequence $(|\mathscr{PT}_m^4|)_{m=0}^{+\infty}$ is surprisingly the sequence A000129 of Pell numbers described in OEIS.

## Acknowledgment

The authors would like to thank the anonymous referees, whose suggestions improved the paper.

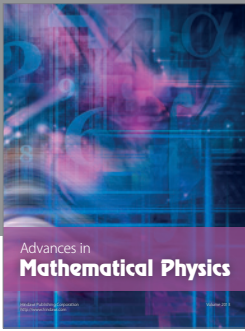 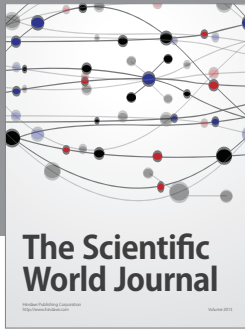 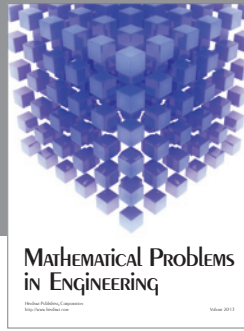 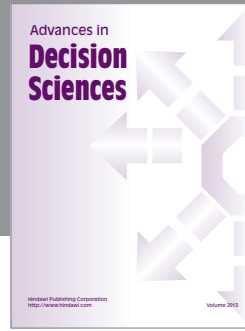 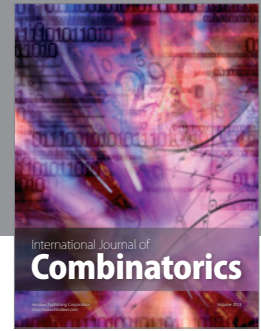
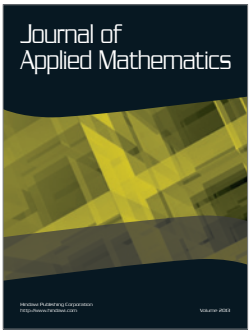 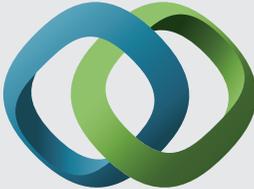 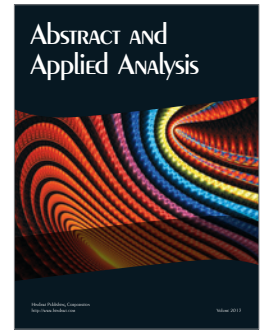
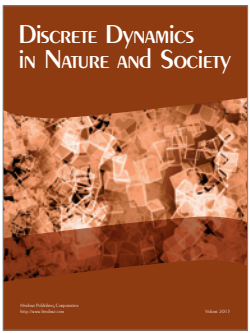 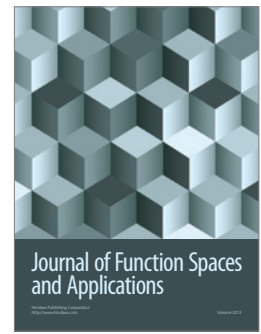
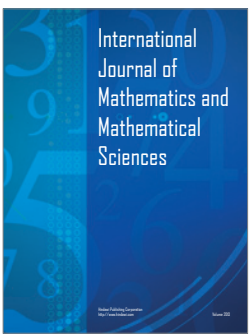 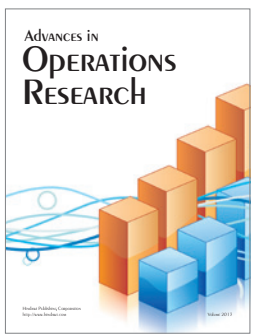 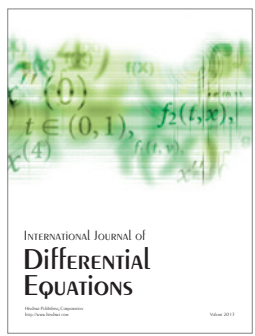 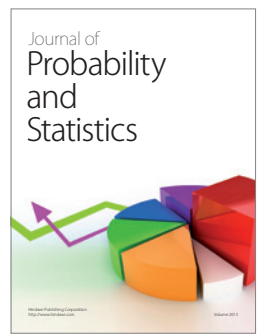 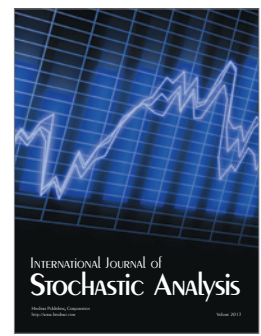
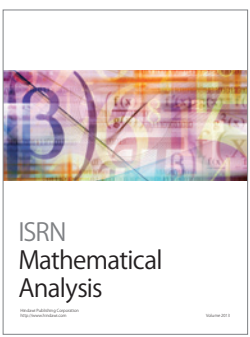 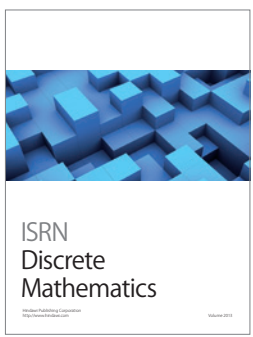 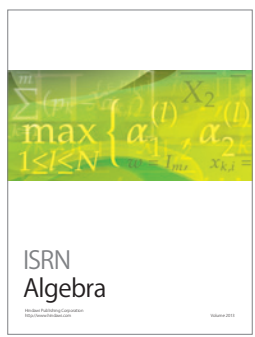 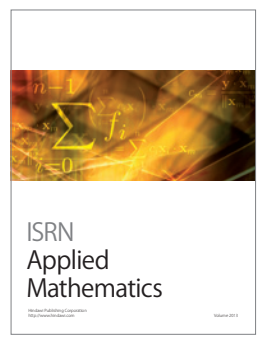 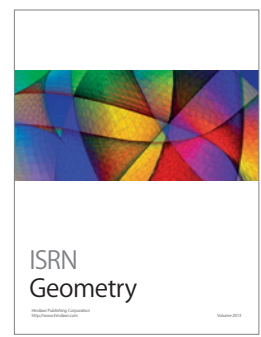

Submit your manuscripts at
http://www.hindawi.com